\documentclass[conference]{IEEEtran}
\IEEEoverridecommandlockouts
\usepackage{cite}
\usepackage{amsmath,amssymb,amsfonts}
\usepackage{algorithmic}
\usepackage{graphicx}
\usepackage{textcomp}
\usepackage{mystyle}

\def\BibTeX{{\rm B\kern-.05em{\sc i\kern-.025em b}\kern-.08em
		T\kern-.1667em\lower.7ex\hbox{E}\kern-.125emX}}
\begin{document}
	
	\title{An Extension of Averaged-Operator-Based Algorithms\\
		%{\footnotesize \textsuperscript{*}Note: Sub-titles are not captured in Xplore and
		%should not be used}	
		\thanks{This work was supported by the Funda\c{c}\~{a}o para a Ci\^{e}ncia e Tecnologia within the Portuguese Ministry for Science, Technology and Higher Education under Project UID/EEA/50008/2013 and Grant BPD/N.º 134 - 16/10/2017.}
	}
	
	\author{\IEEEauthorblockN{Miguel Sim\~{o}es}
		\IEEEauthorblockA{\textit{Instituto de Telecomunica\c{c}\~{o}es} \\
			\textit{Instituto Superior T\'{e}cnico, Univ. Lisboa}\\
			Lisbon, Portugal \\
			miguel.simoes@lx.it.pt}
		\and
		\IEEEauthorblockN{Jos\'{e} Bioucas-Dias}
		\IEEEauthorblockA{\textit{Instituto de Telecomunica\c{c}\~{o}es} \\
			\textit{Instituto Superior T\'{e}cnico, Univ. Lisboa}\\
			Lisbon, Portugal \\
			bioucas@lx.it.pt}
		\and
		\IEEEauthorblockN{Luis B. Almeida}
		\IEEEauthorblockA{\textit{Instituto de Telecomunica\c{c}\~{o}es} \\
		\textit{Instituto Superior T\'{e}cnico, Univ. Lisboa}\\
		Lisbon, Portugal \\
		luis.almeida@lx.it.pt}
		%\and
		%\IEEEauthorblockN{4\textsuperscript{th} Given Name Surname}
		%\IEEEauthorblockA{\textit{dept. name of organization (of Aff.)} \\
		%\textit{name of organization (of Aff.)}\\
		%City, Country \\
		%email address}
		%\and
		%\IEEEauthorblockN{5\textsuperscript{th} Given Name Surname}
		%\IEEEauthorblockA{\textit{dept. name of organization (of Aff.)} \\
		%\textit{name of organization (of Aff.)}\\
		%City, Country \\
		%email address}
		%\and
		%\IEEEauthorblockN{6\textsuperscript{th} Given Name Surname}
		%\IEEEauthorblockA{\textit{dept. name of organization (of Aff.)} \\
		%\textit{name of organization (of Aff.)}\\
		%City, Country \\
		%email address}
	}
	
	\maketitle
	
	\begin{abstract}
		Many of the algorithms used to solve minimization problems with sparsity-inducing regularizers are generic in the sense that they do not take into account the sparsity of the solution in any particular way. However, algorithms known as semismooth Newton are able to take advantage of this sparsity to accelerate their convergence. We show how to extend these algorithms in different directions, and study the convergence of the resulting algorithms by showing that they are a particular case of an extension of the well-known Krasnosel'ski\u{\i}--Mann scheme. %Some alternative interpretations of these methods, as well as some applications, are also discussed. In particular, the methods are experimentally shown to be able to achieve substantially faster convergence rates than standard first- and second-order methods when solving a simple problem.
	\end{abstract}
	
	
	\section{Introduction}
			
	\subsection{Background}
			
	The objective functions of many signal-processing problems can be formulated as sums of two proper lower-semicontinuous convex functions: one that is smooth, $f: \mathbb{R}^n \to ]-\infty,+\infty]$, and another one that need not be smooth, $g: \mathbb{R}^n \to ]-\infty,+\infty]$. The resulting problem is
	\begin{equation} \label{eq:split}
	\underset{\mathbf{x}  \in \mathbb{R}^n}{\text{minimize}} \quad f (\mathbf{x}) + g (\mathbf{x}).
	\end{equation}
	Such problems are typically large-scale and can be solved by using splitting methods, which convert~\eqref{eq:split} into a sequence of separable subproblems. The (relaxed) \emph{forward--backward} method~\cite{Figueiredo2003, Daubechies2004} is an example of such methods. Its iterations can be broken into a gradient (forward) step on $f$ and a proximal (backward) step on $g$, performed consecutively---see Algorithm~\ref{algo:rfba}, where $\text{prox}_{\tau g}$ denotes the \emph{proximal operator} of function $g$, i.e., $\text{prox}_{\tau g}(\mathbf{x}) \eqdef \text{arg min}_{\mathbf{u} \in \mathbb{R}^n} \left\{ g(\mathbf{u}) + \frac{1}{2 \tau} \| \mathbf{x} - \mathbf{u} \|^2 \right\}$~\cite{Moreau1962}.
	\begin{algorithm} \label{algo:rfba}
		Choose $\mathbf{x}^0 \in \mathbb{R}^n,\ \tau > 0$\;
		$k \leftarrow 1$\;
		\While{stopping criterion is not satisfied}{
			Choose $\lambda^k > 0$\;
			$\mathbf{x}^{k + 1} \leftarrow \mathbf{x}^k + \lambda^k \left( \text{prox}_{\tau g} \left(\mathbf{x}^k - \tau \nabla f \left(\mathbf{x}^k \right)\right) - \mathbf{x}^k \right)$\; \label{eq:iter_rfba}
			$k \leftarrow k + 1$\;
		}
		\caption{Relaxed forward--backward method.}
	\end{algorithm}

	When analyzing the properties of many of these and other algorithms, it can be advantageous to use the theory of monotone operators~\cite{Byrne2004}. Let $2^{\mathbb{R}^n}$ denote the {power} set of $\mathbb{R}^n$. A set-valued operator $A: \mathbb{R}^n \to 2^{\mathbb{R}^n}$ is said to be \emph{monotone} if $\langle \mathbf{u}-\mathbf{v}, \mathbf{x}-\mathbf{y} \rangle \geq 0$ for all $(\mathbf{x},\mathbf{u}) \in \text{gra } A$ and $(\mathbf{y},\mathbf{v}) \in \text{gra } A$, where $\text{gra } A$ denotes the {graph} of A, and it is said to be \emph{maximally monotone} if there exists no other monotone operator whose graph properly contains $\text{gra } A$. Monotone operators are connected to optimization problems as follows. Take, for example, \eqref{eq:split}. According to Fermat's rule, its solutions should satisfy the inclusion $0 \in \nabla f (\mathbf{x}) + \partial g (\mathbf{x})$, where the set-valued operator $\partial g : \mathbb{R}^n \to 2^{\mathbb{R}^n} : \mathbf{x} \to \partial g (\mathbf{x})$ denotes the \emph{subdifferential} of $g$ (in the sense of Moreau and Rockafellar~\cite[Chapter 23]{Rockafellar1970}). The operators $\nabla f$ and $\partial g$ are examples of maximally-monotone operators~\cite[Theorem 20.40]{Bauschke2011}. Problem~\eqref{eq:split} can be seen as a particular case of the problem of finding a zero of the sum of two monotone operators $A$ and $C$, i.e.,
	\begin{equation} \label{eq:inclusion}
	\text{find } \mathbf{x} \in \mathbb{R}^n \quad \text{such that } 0 \in A \left(\mathbf{x}\right) + C \left(\mathbf{x}\right),
	\end{equation}
	if one makes $A=\partial g$ and $C=\nabla f$. Problem~\eqref{eq:inclusion} may be solved using a generalized version of Algorithm~\ref{algo:rfba}, in which Line~\ref{eq:iter_rfba} is replaced with
	\begin{equation} \label{algo:rfba_op}
		\mathbf{x}^{k + 1} \leftarrow \mathbf{x}^k + \lambda^k \left( J_{\tau A} \left( \mathbf{x}^k - \tau C \left(\mathbf{x}^k\right) \right) - \mathbf{x}^k \right),
	\end{equation}
	where $J_{\tau A} \eqdef (\text{Id} + \tau A)^{-1}$ is the \emph{resolvent} of operator $A$ and $\text{Id}$ denotes the \emph{identity} operator. Note that $J_{\tau \partial g} = \text{prox}_{\tau g}$ \cite[Example 23.3]{Bauschke2011}. 
	
	Problem~\eqref{eq:inclusion} can alternatively be written as the problem of finding a fixed point of the operator $R \eqdef J_{\tau A} \circ (\text{Id} - \tau C)$:
	\begin{equation} \label{eq:mon_incl_T}
	\text{find } \varx \in \mathbb{R}^n \quad \text{such that } \neop \left( \varx \right) = \varx.
	\end{equation}
	In general, the solutions of a convex optimization problem correspond to the fixed points of a certain operator, and an iterative optimization algorithm corresponds to a fixed-point method. %There are at least two advantages in looking at convex problems and at these algorithms as particular cases of monotone inclusions and of fixed-point methods, respectively: (a) the analysis of the convergence of these algorithms becomes easier, since one is able to take advantage of the many results on operator theory that are available from the literature, and (b) besides convex problems, other problems of interest can also be seen as particular cases of monotone inclusions. 
	We can rewrite~\eqref{algo:rfba_op} as
	\begin{equation} \label{algo:km}
		\mathbf{x}^{k + 1} \leftarrow T_{\lambda^k} \left( \mathbf{x}^k \right) \eqdef \mathbf{x}^k + \lambda^k (R \left( \mathbf{x}^k \right) - \mathbf{x}^k).
	\end{equation}
	We say that an operator $R: \mathbb{R}^n \to \mathbb{R}^n$ is nonexpansive if $\| \mathbf{u} - \mathbf{v} \| \leq \| \mathbf{x} - \mathbf{y} \|$ for all $(\mathbf{x},\mathbf{u}) \in \text{gra } R$ and $(\mathbf{y},\mathbf{v}) \in \text{gra } R$. Let $R$ be a generic nonexpansive operator and let $\lambda \in \; ]0,1[$. Then the operator $T \eqdef \left(\eye - \lambda \right) + \lambda R$ is said to be $\lambda$-\emph{averaged}. It obeys the following contractive property~\cite[Proposition 4.25]{Bauschke2011}:
	\begin{align}
		&\norm{T \left(\varx\right) - T \left(\vary\right)}{}^2 \nonumber \\
		& \, \leq \norm{\varx - \vary}{}^2 - \frac{1 - \lambda}{\lambda} \norm{\left( \text{Id} - T \right) \left(\varx\right) - \left( \text{Id} - T \right) \left(\vary\right)}{}^2 \label{eq:contract_avop}
	\end{align}
	for all $\varx, \, \vary \in \mathbb{R}^n$. In particular, when $\lambda = 1/2$, $T$ is said to be \emph{firmly nonexpansive}. The resolvents of maximally-monotone operators are firmly-nonexpansive~\cite[Corollary 23.8]{Bauschke2011}. Iteration~\eqref{algo:km} is known as the \emph{Krasnosel'ski\u{\i}--Mann} scheme and is the basis of not only the forward--backward method but also other optimization algorithms, such as the Douglas--Rachford one~\cite{Byrne2004, Bauschke2011}. It can be shown that, under certain conditions, the Krasnosel'ski\u{\i}--Mann scheme converges to $\text{Fix } R$, where $\text{Fix } R$ denotes the set of fixed points of $R$.
	 
	The convergence rate of the forward--backward method (Algorithm~\ref{algo:rfba}) can be shown to be sublinear, or, under certain assumptions, to be linear. This rate can often be improved by incorporating \emph{second-order} information about $f$ if this function is twice-differentiable. The local convergence rate of second-order methods is superlinear or even quadratic. As an example, consider the second-order version of Algorithm~\ref{algo:rfba}, which is given by replacing Line~\ref{eq:iter_rfba} with the iteration $\mathbf{x}^{k + 1} \leftarrow \mathbf{x}^k + \lambda^k \left( \text{prox}_g^{\mathbf{B}^k} \left( \mathbf{x}^k - \left[ \mathbf{B}^k \right]^{-1} \nabla f \left( \mathbf{x}^k \right) \right) - \mathbf{x}^k \right)$~\cite{Schmidt2011, Becker2012, Lee2014}, where $\mathbf{B}^k$ is a \ac{PD} matrix (the Hessian of $f$ or an approximation of it) and $\text{prox}_g^{\mathbf{B}^k}$ denotes the proximal operator of $g$ relative to the norm $\| \cdot \|^2_{\mathbf{B}^k}$, i.e., $\text{prox}_{g}^{\mathbf{B}^k} (\mathbf{x}) \eqdef \text{arg min}_{\mathbf{u} \in \mathbb{R}^n} \left\{ g(\mathbf{u}) + \frac{1}{2} \| \mathbf{x} - \mathbf{u} \|^2_{\mathbf{B}^k} \right\}$. More generally, and from an operator-centric perspective, by using second-order methods such as these, one is actually solving a left-preconditioned version of~\eqref{eq:inclusion}, in the sense that instead of directly tackling that problem we are considering problems that share the same set of solutions but may be more convenient to solve:
	\begin{equation} \label{eq:mon_incl_precond}
	\text{find } \varx \in \mathbb{R}^n \quad \text{such that } 0 \in \fbvm \maxmonresol \left( \varx \right) + \fbvm \cocoercop \left( \varx \right),
	\end{equation}
	where $\fbvm$ is a \ac{PD} operator. In what follows, we denote positive definiteness by $\fbvm \succ 0$ and positive semidefiniteness by $\fbvm \succeq 0$.
	
	\subsection{Contributions}
	
	The basis of this work is the study of the following alternative scheme to~\eqref{algo:km}:
	\begin{equation} \label{eq:fixed_point_iter_opwavop}
	\varx^\iite = \avop{\ave^\ite} \left( \varx^\ite \right) \eqdef \varx^\ite + \ave^\ite \left( \neop \left( \varx^\ite \right) - \varx^\ite \right),
	\end{equation}
	where, for every $\ite$, $\ave^\ite$ is a linear operator such that $\eye \succ \ave^\ite \succ 0$. For convenience, we call the operators $\avop{\ave^\ite}$, \emph{operator-weighted averaged operators}. It is clear that if, for all $k$, we make $\ave^\ite = \lambda^\ite \eye$, we recover~\eqref{algo:km}. 
	
	Iteration~\eqref{eq:fixed_point_iter_opwavop} can be interpreted in different ways. For example, if $\ave^\ite$ is fixed, i.e., if, for all $\ite$, $\ave^\ite = \ave$, where $\ave \succ 0$, that iteration can also be seen as a left-preconditioning scheme to solve~\eqref{eq:mon_incl_T}: 
	\begin{equation} \label{eq:mon_incl_precond_T}
	\text{find } \varx \in \mathbb{R}^n \quad \text{such that } \ave \neop \left( \varx \right) = \ave \varx.
	\end{equation}
	
	\subsection{Notation and outline}
	
	A detailed account of the notions listed in this section can be found in the work of Bauschke and Combettes~\cite{Bauschke2011}. We denote the \emph{scalar product} of a Hilbert space by $\langle \cdot , \cdot \rangle$ and the associated \emph{norm} by $\| \cdot \|$. The \emph{range} of an operator $A$ is denoted by $\text{ran } A$, and the \emph{adjoint} of $A$ by $A^*$. We say that an operator $A : \mathbb{R}^n \to \mathbb{R}^n$ is \emph{Lipschitz continuous} with constant $L > 0$ if $\| \mathbf{u} - \mathbf{v} \| \leq L \| \mathbf{x} - \mathbf{y} \|$, for all $(\mathbf{x},\mathbf{u}) \in \text{gra } A$ and $(\mathbf{y},\mathbf{v}) \in \text{gra } A$. Additionally, let $\Gamma_0 (\mathbb{R}^n)$ denote the class of all proper lower-semicontinuous convex functions from $\mathbb{R}^n$ to $]{-\infty},+\infty]$. Given two functions $f \in \Gamma_0 (\mathbb{R}^n)$ and $g \in \Gamma_0 (\mathbb{R}^n)$, their \emph{infimal convolution} is denoted by  $\infconv{f}{g}$. The Legendre--Fenchel \emph{conjugate} of a function $f$ is denoted by $f^*$. The \emph{indicator function} of a set $C \in \mathbb{R}^{n}$ is defined as $\delta_C(\mathbf{x}) \eqdef 0$ if $\mathbf{x} \in C$, $\delta_C(\mathbf{x}) \eqdef + \infty$ otherwise. We use the notation $\{\mathbf{x}^k\}$ as a shorthand for representing the sequence $\{\mathbf{x}^k\}_{k=1}^{+\infty}$. The space of \emph{absolutely-summable sequences} in $\mathbb{R}$ is denoted by $\spaceseq^1 (\mathbb{N})$; the set of summable sequences in $[0, + \infty [$ is denoted by $\spaceseq^1_+ (\mathbb{N})$. Bold lowercase letters denote vectors and bold uppercase letters denote matrices. $[\mathbf{a}]_i$ denotes the $i$-th element of a vector $\mathbf{a}$, $[\mathbf{A}]_{:j}$ denotes the $j$-th column of a matrix $\mathbf{A}$, and $[\mathbf{A}]_{ij}$ denotes the element in the $i$-th row and $j$-th column of a matrix $\mathbf{A}$. $\mathbf 0$ denotes a \emph{zero} vector or matrix of appropriate size. The maximum and signum operators are denoted by $\text{max}(\cdot)$ and $\text{sgn}(\cdot)$, respectively.
	
	%Although the idea behind operator-weighted averaged operators is quite natural, it has not, to the best of our knowledge, been formalized before. The asymptotic analysis of the fixed-point iterations of these operators seems to be new; it is based on some previous results on the convergence of variable-metric quasi-Fej\'{e}r sequences~\cite{Combettes2013,Combettes2014b}. 
	
	The structure of this work is as follows. In Section~\ref{sec:ssn}, we briefly discuss a class of algorithms known as semismooth Newton methods. In Section~\ref{sec:extension}, we study the scheme given by~\eqref{eq:fixed_point_iter_opwavop}, and show how it can be used to solve a primal--dual problem first studied by Combettes and Pesquet~\cite{Combettes2011b}. In Section~\ref{sec:apps}, we present a simple application of the proposed method to solve an inverse problem. Section~\ref{sec:conclusions} concludes. Due to space constraints, we omit the proofs of the results discussed in Section~\ref{sec:extension}; these proofs can be consulted elsewhere~\cite[Chapter 5]{Simoes2017}.
	
	\section{Semismooth Newton methods} \label{sec:ssn}
	
	\emph{Semismooth Newton} methods were originally developed with the goal of using Newton-like methods to minimize certain nonsmooth functions at a superlinear convergence rate. To illustrate why these methods may be useful when solving problems of the form of~\eqref{eq:split}, consider, as an example, that $f = \| \mathbf{y} - \mathbf{H} \cdot \|^2$, and $g = \mu \| \cdot \|_1$, where $\mathbf{y} \in \mathbb{R}^m$, $\mathbf{H} \in \mathbb{R}^{m \times n}$, and $\mu > 0$. For problems such as these, it was shown by Hinterm\"{u}ller~\cite{Hintermuller2003} that some semismooth Newton methods are equivalent to some active-set methods. As we discuss in Section~\ref{sec:apps}, the fact that these methods can be written as active-set ones allows for significant time savings when solving certain problems, namely the ones involving sparsity-inducing regularizers, as is the case of the $\ell_1$ norm. 
	
	Let $G: \mathbb{R}^n \to \mathbb{R}^n$ be an operator such that $G: \mathbf{x} \to \mathbf{x} - \text{prox}_{\mu \| \mathbf{x} \|_1} \left(\mathbf{x} - 2 \mu \mathbf{H}^* (\mathbf{H} \mathbf{x} - \mathbf{y}) \right)$. The solution of the problem under consideration should satisfy the nonlinear equation $G ({\mathbf{x}}) = \mathbf{0}$, which is nonsmooth, since $\text{prox}_{\mu \| \cdot \|_1}$ is not everywhere differentiable. There are, however, generalizations of the concept of differentiability that are applicable to an operator such as $G$. One of them is the B(ouligand)-differential~\cite[Definition 4.6.2]{Facchinei2003}, which is defined as follows. Suppose that a generic operator $G: D \subset \mathbb{R}^n \to \mathbb{R}^m$ is locally Lipschitz, where $D$ is an open subset. Then by Rademacher's theorem, $G$ is differentiable almost everywhere in $D$. Let $C$ denote the subset of $\mathbb{R}^n$ consisting of the points where $G$ is differentiable (in the sense of Fr\'{e}chet~\cite[Definition 2.45]{Bauschke2011}). The B-differential of $G$ at $\mathbf{x}$ is $\partial_B \, G (\mathbf{x}) \eqdef \left\{ \lim_{\mathbf{x}^j \to \mathbf{x}} \nabla G \left(\mathbf{x}^j \right) \right\}$, where $\{\mathbf{x}^j\}$ is a sequence such that $\mathbf{x}^j \in C$ for all $j$ and $\nabla G(\mathbf{x}^j)$ denotes the Jacobian of $G$ at $\mathbf{x}^j$. 
	
	The B-differential of an operator at a given point may not be unique: for example, take $\text{prox}_{\mu \| \cdot \|_1} (\mathbf{x})$, which can be evaluated element-wise by computing $\text{max} \left\{ \big| [\mathbf{x}]_i \big| - \mu, 0 \right\} \circ \text{sgn} \left([\mathbf{x}]_i \right)$ for $i \in \{1, \cdots, n\}$. A possible $\mathbf{H} \in \partial_B \, \text{prox}_{\mu \| \cdot \|_1} (\mathbf{x})$ is a binary diagonal matrix defined as~\cite[Proposition 3.3]{Griesse2008}
	\begin{equation} \label{eq:partial_B_G}
		[\mathbf{H}]_{ii} = \begin{cases}
		1 & \text{if } \big| [\mathbf{x}]_i \big| > \mu,\\
		0 & \text{otherwise.}
		\end{cases}
	\end{equation}
	
	This generalization of the concept of differentiability can also be used to formulate the so-called semismooth Newton method, which is characterized by the iteration $\mathbf{x}^{k+1} \leftarrow \mathbf{x}^{k} - [\mathbf{H}^k]^{-1} \, G(\mathbf{x}^k)$, where $\mathbf{H}^k \in \partial_B  \, G (\mathbf{x}^k)$. It can be shown that this method locally converges superlinearly for operators known as semismooth~\cite{Qi1993b}. Let $\mathbf{x} \in D$ and $\mathbf{d} \in \mathbb{R}^n$; semismooth operators are operators that are directionally differentiable at a neighborhood of $\mathbf{x}$ and that, for any $\mathbf{H} \in \partial_B \, G (\mathbf{x}+\mathbf{d})$, satisfy the condition $\mathbf{H} \mathbf{d} - G'(\mathbf{x}; \mathbf{d})=o(\|\mathbf{d}\|)$ for $\mathbf{d} \to \mathbf{0}$, where $G'(\mathbf{x}; \mathbf{d})$ denotes the directional derivative of $G$ at $\mathbf{x}$ along $\mathbf{d}$. Examples of semismooth functions are the Euclidean norm and piecewise-differentiable functions~\cite[Chapter 2]{Ulbrich2011}, $\text{prox}_{\mu \| \cdot \|_1} (\mathbf{x})$ being an example of the latter. Note that the semismooth Newton method is a particular case of~\eqref{eq:fixed_point_iter_opwavop}, although we impose that $\eye \succ \ave^\ite \succ 0$ in the latter equation, which is not necessarily true for this method.

	\section{An extension of averaged-operator-based algorithms} \label{sec:extension}

	In this section, we define operator-weighted averaged operators, and show that they have a contractive property. We also study the asymptotic behavior of fixed-point iterations of these operators. Such iterations can be seen as an extension of the Krasnosel'ski\u{\i}--Mann scheme [cf.~\eqref{algo:km}]. We base our analysis on the fact that these iterations produce a sequence that is variable-metric Fej\'{e}r monotone~\cite{Combettes2013,Combettes2014b}. We then present an algorithm that uses operator-weighted averaged operators, and that solves a primal--dual problem that encapsulates many problem formulations~\cite{Combettes2011b, Combettes2014b}.

	\subsection{An extension of the Krasnosel'ski\u{\i}--Mann scheme} \label{sec:eKM}
	
	\begin{definition}[Operator-weighted averaged operators]
		Let $\nesubset$ be a nonempty subset of $\mathbb{R}^n$, let $\eps \in \; ]0, 1[$, and let $\ave$ be an operator in $\mathbb{R}^n$ such that 
		\begin{equation} \label{eq:ave_ass}
		\Meig \eye \succeq \ave \succeq \meig \eye, \quad \text{where } \Meig, \, \meig \in [\eps, 1-\eps].
		\end{equation}
		We say that an operator $\avop{\ave}: \nesubset \to \mathbb{R}^n$ is an operator-weighted averaged operator if there exists a nonexpansive operator $\neop: \nesubset \to \mathbb{R}^n$ such that 
		\begin{equation} \label{eq:aveop}
		\avop{\ave} \define (\eye - \ave) + \ave \neop.
		\end{equation}
	\end{definition}
	
	We have proved the following results:
	
	\begin{proposition} \label{th:contractive}
		Let $\nesubset$ be a nonempty subset of $\mathbb{R}^n$, let $\eps \in \; ]0, 1[$, let $\ave$ be an operator in $\mathbb{R}^n$ satisfying~\eqref{eq:ave_ass}, let $\neop: \nesubset \to \mathbb{R}^n$ be a nonexpansive operator, and let $\avop{\ave}: \nesubset \to \mathbb{R}^n$ be an operator as defined in~\eqref{eq:aveop}. Then the operator $\avop{\ave}$ is $\Meig$-averaged in the metric induced by $\inv{\ave}$. In other words, the operator $\avop{\ave}$ verifies
		\begin{align*}
		&\norm{\avop{\ave} \left( \varx \right) - \avop{\ave} \left( \vary \right)}{\inv{\ave}}^2 \\ \nonumber
		&\leq \norm{\varx - \vary}{\inv{\ave}}^2 - \frac{1 - \Meig}{\Meig} \norm{\left( \emph{\eye} - \avop{\ave} \right) \left( \varx \right) - \left( \emph{\eye} - \avop{\ave} \right) \left( \vary \right)}{\inv{\ave}}^2
		\end{align*}
		for all $\varx, \, \vary \in \nesubset$.
	\end{proposition}
	
%	Consider iterarions of operator-weighed averaged operators of the form of iteration~\eqref{eq:fixed_point_iter_opwavop}. The following theorem establishes some convergence properties of this algorithm.
%	
	\begin{theorem} \label{th:eKM}
		Let $\nesubset$ be a nonempty closed convex subset of $\mathbb{R}^n$, let $\eps \in \; ]0, 1[$, let $\seq{\seqave^\ite} \in \spaceseq^1_+(\mathbb{N})$, let $\seq{\ave^\ite}$ be a sequence of \ac{PD} operators in $\mathbb{R}^{n \times n}$ such that, for all $\ite \in \mathbb{N}$, 
		\begin{equation} \label{eq:seq_ave_ass}
		\begin{cases}
		\Meig^\ite \emph{\eye} \succeq \ave^\ite \succeq \meig^\ite \emph{\eye},\\
		\Meig^\ite, \, \meig^\ite \in [\eps, 1-\eps],\\
		\left( 1 + \seqave^\ite \right) \ave^\iite \succeq \ave^\ite,
		\end{cases}
		\end{equation}
		and let $\neop: \nesubset \to \nesubset$ be a nonexpansive operator such that $\emph{\fix} \neop \neq \emptyset$. Additionally, let $\varx^0 \in \nesubset$ and let, for all $\ite$, $\seq{\varx^\ite}$ be a sequence generated by~\eqref{eq:fixed_point_iter_opwavop}. Then $\seq{\varx^\ite}$ converges to a point in $\emph{\fix} \neop$.
%		\end{enumerate}
	\end{theorem}
	
%	\begin{proof}
%		We defer all proofs to Section~\ref{sec:proofs}.
%	\end{proof}
	
	\subsection{Primal--dual optimization algorithms} \label{sec:pdprob}
	
	Combettes and Pesquet studied a primal--dual problem that generalizes many problems~\cite[Problem 4.1]{Combettes2011b}. By being able to devise an algorithm to solve this problem, we are effectively tackling a large number of problems simultaneously (problem~\eqref{eq:split} is one of these). Let $m$, $n$, and $\dimr$ be strictly-positive integers, let $\cvxone \in \lsc(\mathbb{R}^n)$, let $\coco \in \; ]0, + \infty[$, let $\smooth: \mathbb{R}^n \to ]{-\infty},+\infty]$ be convex and differentiable with a $\inv{\coco}$-Lipschitzian gradient, and let $\biasvar \in \mathbb{R}^n$. For every $\itr \in \{ 1, \dots, \dimr \}$, let $\biasvarlnop_\itr \in \mathbb{R}^{m_\itr}$, let $\cvxtwo_\itr \in \lsc(\mathbb{R}^{m_\itr})$, let $\strmaxmonopparam_\itr \in \; ]0, +\infty[$, let $\strong_\itr \in \lsc(\mathbb{R}^{m_\itr})$ be $\strmaxmonopparam_\itr$-strongly convex,\footnote{A function $\strong$ is said to be $\strmaxmonopparam$-\emph{strongly convex} if $\strong - \frac{\strmaxmonopparam}{2} \langle \mathbf{x}, \mathbf{x} \rangle$ is convex, for some $\strmaxmonopparam > 0$.} let $\lnop_\itr \in \mathbb{R}^{m_\itr \times n}$ such that $\lnop_\itr \neq 0$, and let $\pdomega_\itr$ be real numbers in $]0, 1]$ such that $\sum_{\itr = 1}^{\dimr} \pdomega_\itr = 1$. The problem is as follows:
	
	\begin{problem} \label{th:problem}

		Solve the primal minimization problem,
		\begin{equation*} \label{eq:primalproblem2}
			\underset{\pr \in \mathbb{R}^n}{\text{{minimize}}} \, \cvxone (\pr) + \sum_{\itr=1}^{\dimr} \pdomega_\itr \left(\infconv{\cvxtwo_\itr}{\strong_\itr} \right) \left(\lnop_\itr \pr - \biasvarlnop_\itr \right) + \smooth (\pr) - \innerpro{\pr}{\biasvar}{},
		\end{equation*}
		together with its corresponding dual minimization problem,
		\begin{align*} 
			\underset{\du_1 \in \mathbb{R}^{m_1}, \cdots, \du_\itr \in \mathbb{R}^{m_\itr}}{\text{{minimize}}} &\, \left( \infconv{\conj{\cvxone}}{\conj{\cvxtwo}} \right) \left( \biasvar - \sum_{\itr=1}^{\dimr} \pdomega_\itr  \conj{\lnop_\itr} \du_\itr \right) \\
			&\, + \sum_{\itr=1}^{\dimr} \pdomega_\itr \left( \conj{\cvxtwo_\itr} (\du_\itr) + \conj{\strong_\itr} (\du_\itr) + \innerpro{\du_\itr}{\biasvarlnop_\itr}{} \right). \label{eq:dualproblem2}
		\end{align*}
		The sets of solutions to these primal and dual problems are denoted by $P$ and $D$, respectively.
		
	\end{problem}
	
	Consider Algorithm~\ref{algo:stackevmfbapp} to solve Problem~\ref{th:problem}. In what follows, for all $\itr$, $\seq{\vm^\ite}$, $\seq{\ave^\ite}$, $\seq{\vm^\ite_\itr}$, $\seq{\ave^\ite_\itr}$ are sequences of linear operators, and $\seq{\errorpr^\ite}$, $\seq{\errordu^\ite_\itr}$, $\seq{\errorresolpr^\ite}$, $\seq{\errorresoldu^\ite_\itr}$ are absolutely-summable sequences that can be used to model errors. Algorithm~\ref{algo:stackevmfbapp} is an extension of~\cite[Example 6.4]{Combettes2014b}.
	\begin{algorithm} \label{algo:stackevmfbapp}
		Choose $\pr^0 \in \mathbb{R}^n$ and $\du_1^0 \in \mathbb{R}^{m_1}, \cdots, \du_\itr^0 \in \mathbb{R}^{m_\itr}$\;
		$\ite \leftarrow 1$\;		
		\While{stopping criterion is not satisfied}{
			\For{$\itr = 1, \dots, \dimr$}{
				Choose $\vm^\ite_\itr, \, \ave^\ite_\itr \succ 0 \text{ s.t. } \ave^\ite_\itr \prec \eye$\;
				$\pdu^\ite_\itr =\prox^{\inv{(\vm_\itr^\ite)}}_{\conj{\cvxtwo}_\itr} \big( \du^\ite_\itr + \vm_\itr^\ite \big( \lnop_\itr \pr^\ite$
					\newline\makebox[3cm]{}$ - \grad{\conj{\strong}_\itr} \left(\du^\ite\right) - \errorresoldu^\ite_\itr - \biasvarlnop_\itr \big) \big) + \errordu^\ite_\itr$\;
				$\duo^\ite_\itr = 2 \pdu^\ite_\itr - \du^\ite_\itr$\;
				$\du^\iite_\itr = \du^\ite_\itr + \ave_\itr^\ite \left( \pdu^\ite_\itr - \du^\ite_\itr \right)$\;
			}
			Choose $\vm^\ite, \, \ave^\ite \succ 0 \text{ s.t. } \ave^\ite \prec \eye$\;
			$\ppr^\ite =$
				$\prox_{\cvxone}^{\inv{(\vm^\ite)}} \big( \pr^\ite - \vm^\ite \big( \sum_{\itr = 1}^{\dimr} \pdomega_\itr \conj{\lnop}_\itr \duo^\ite_\itr$
				\newline\makebox[3cm]{}$+ \grad{\smooth} \left(\pr^\ite\right) + \errorresolpr^\ite - \biasvar \big) \big) + \errorpr^\ite$\;
			$\pr^\iite =  \pr^\ite + \ave^\ite \left(\ppr^\ite - \pr^\ite \right)$\; \label{eq:ssn_in_algo}
			$\ite \leftarrow \iite$\;		
		}
		\caption{An application of~\eqref{eq:fixed_point_iter_opwavop} to solve Problem~\ref{th:problem}.}
	\end{algorithm}
	The following corollary establishes some convergence properties of Algorithm~\ref{algo:stackevmfbapp}.
	
	\begin{corollary} \label{th:stackevmfbapp}
		Suppose that 
		\begin{equation*} \label{eq:assumption2}
		\biasvar \in \emph{\ran} \left( \subgrad{\cvxone} + \sum_{\itr=1}^{\dimr} \pdomega_\itr \conj{\lnop_\itr} \left( \infconv{\subgrad{\cvxtwo}_\itr}{\subgrad{\strong}_\itr} \right) \left(\lnop_\itr \cdot - \biasvarlnop_\itr \right) + \grad{\smooth} \right)
		\end{equation*}	
		and set $\mincocomon \eqdef \min \{\coco, \strmaxmonopparam_1, \dots, \strmaxmonopparam_\dimr \}$. Let $\seq{\vm^\ite}$ be a sequence of \ac{PD} operators in $\mathbb{R}^{n \times n}$ and, for every $\itr \in \{ 1, \dots, \dimr \}$, let $\seq{\vm^\ite_\itr}$ be a sequence of \ac{PD} operators in $\mathbb{R}^{m_\itr \times m_\itr}$ such that, for all $\ite \in \mathbb{N}$,
		\begin{equation} \label{eq:stackevmfbass1}
		\begin{cases}
		\Meigvm \emph{\eye} \succeq \vm^\ite \succeq \meigvm \emph{\eye},\\
		\Meigvm \emph{\eye} \succeq \vm^\ite_\itr \succeq \meigvm \emph{\eye},\\
		\Meigvm, \, \meigvm \in \ ]0, + \infty[,
		\end{cases}
		\end{equation}
		let $\fbeps \in \;]0, \min \{ 1, \mincocomon \}[$, let $\seq{\ave^\ite}$ be a sequence of \ac{PD} operators in $\mathbb{R}^{n \times n}$, and let $\seq{\ave^\ite_\itr}$ be a sequence of \ac{PD} operators in $\mathbb{R}^{m_\itr \times m_\itr}$ such that, for all $\ite$,
		\begin{equation} \label{eq:stackevmfbassave}
		\begin{cases}
		\ave^\ite \vm^\ite = \vm^\ite \ave^\ite,\\
		\ave^\ite_\itr \vm^\ite_\itr = \vm^\ite_\itr \ave^\ite_\itr,\\
		\Meig \emph{\eye} \succeq \ave^\ite \succeq \meig \emph{\eye},\\
		\Meig \emph{\eye} \succeq \ave^\ite_\itr \succeq \meig \emph{\eye},\\
		\Meig, \, \meig \in [\fbeps, 1],
		\end{cases}
		\text{and }
		\begin{cases}
		\ave^\iite \vm^\iite \succeq \ave^\ite \vm^\ite, \\
		\ave^\iite_\itr \vm^\iite_\itr \succeq \ave^\ite_\itr \vm^\ite_\itr.		
		\end{cases}
		\end{equation} 
			
		Let, for all $\itr$, $\seq{\errorpr^\ite}$, $\seq{\errordu^\ite}$, $\seq{\errorresolpr^\ite_\itr}$, $\seq{\errorresoldu^\ite_\itr} \in \spaceseq^1(\mathbb{N})$. For every $\ite$, set $\pddelta^\ite \define \isquareroot{\sum_{\itr = 1}^{\dimr} \pdomega_\itr  \norm{\squareroot{\vm^\ite_\itr} \lnop_\itr \squareroot{\vm^\ite}}{}^2} - 1$ and suppose that $\pdxi^\ite \define \frac{\pddelta^\ite}{(1+\pddelta^\ite) \Meigvm} \geq \frac{1}{2 \mincocomon - \fbeps}$.
				
		Let $\seq{\pr^\ite}$ be a sequence generated by Algorithm~\ref{algo:stackevmfbapp}. Then $\pr^\ite$ converges to a point in $P$ and $\left(\du^\ite_1, \dots, \du^\ite_\dimr \right)$ converges to a point in $D$.
	\end{corollary}
	
	\section{Experiment} \label{sec:apps}
	
	In this section, we give a practical example of a simple problem that can be solved via Algorithm~\ref{algo:stackevmfbapp}. Consider the constrained problem
	\begin{equation} \label{eq:lasso_constraint}
	\underset{\mathbf{x} \in [c,d]^n}{\text{minimize}} \quad \| \mathbf{b} - \mathbf{H} \mathbf{x} \|^2 + \mu \| \mathbf{x} \|_1,
	\end{equation}
	where $\mathbf{b} \in \mathbb{R}^n$, $c \in \mathbb{R}$, $d \in \mathbb{R}$, $\mu>0$, $\mathbf{H} = {{}^{1}\!/_{n}} \widehat{\mathbf{H}}$, and $\widehat{\mathbf{H}} \in \mathbb{R}^{n \times n}$ is a lower-triangular matrix of ones. Griesse and Lorenz studied a non-constrained, and therefore simpler, version of this problem in the context of inverse integration~\cite[Section 4.1]{Griesse2008}. Problem~\eqref{eq:lasso_constraint} can be solved via Algorithm~\ref{algo:stackevmfbapp} if we let $\gamma > 0$, $\tau >0$ and make $m = n$, $\dimr = 1$, $\lnop_1 = \eye$, $\biasvarlnop_1 = \myvec{0}$, $\biasvar = \myvec{0}$, and, for all $\ite$, $\vm^\ite_1 = \gamma \eye$, $\vm^\ite = \tau \eye$, $\errorresoldu^\ite_1 = \myvec{0}$, $\errordu^\ite_1 = \myvec{0}$, $\ave_{1}^\ite = \eye$, $\errorresolpr^\ite = \myvec{0}$, $\errorpr^\ite = \myvec{0}$, $\smooth = \| \mathbf{b} - \mathbf{H} \cdot \|^2$, $\cvxone=\mu \| \cdot \|_1$, $\cvxtwo= \ind{[c,d]^n}{\cdot}$, $\strong_1: \mathbf{u} \to 0$ if $\mathbf{u} = 0$, $\strong_1: \mathbf{u} \to + \infty$ otherwise. 
	
	If we take $\ave^\ite$ to be a sequence of scalars, we recover a version of~\cite[Example 6.4]{Combettes2014b}. However, inspired by the fast convergence properties of the methods discussed in Section~\ref{sec:ssn} and following a similar reasoning to~\cite[Proposition 3.7]{Griesse2008}, we consider the B-differential for the operator $\text{prox}_{\mu \| \cdot \|_1}$ given in~\eqref{eq:partial_B_G} and take $\ave^\ite$ to be the inverse of
	\begin{equation*} \label{eq:op_ssn}
	\inv{\left(\mathbf{P}^\ite\right)} \begin{bmatrix} 
	\tau [\mathbf{H}]^*_{:{I}^k} [\mathbf{H}]^{}_{:{I}^k} & \tau [\mathbf{H}]^*_{:{I}^k} [\mathbf{H}]^{}_{:{A}^k} \\
	\mathbf{0} & \eye
	\end{bmatrix} \mathbf{P}^\ite,
	\end{equation*}
	where
	\begin{equation*}
	\begin{aligned}
	{A}^k &\eqdef \{ i \in \mathbb{N} : \big| \left[\pr^\ite - 2 \tau \left( \mathbf{H}^* \left(\mathbf{H} \mathbf{x}^k - \mathbf{b}\right) + \duo^\ite_1 \right)\right]_i \big| \leq \tau \mu \},\\
	{I}^k &\eqdef \{ i \in \mathbb{N} : \big| \left[\pr^\ite - 2 \tau \left( \mathbf{H}^* \left(\mathbf{H} \mathbf{x}^k - \mathbf{b}\right) + \duo^\ite_1 \right)\right]_i \big| > \tau \mu \},
	\end{aligned}
	\end{equation*}
	and $\seq{\mathbf{P}^\ite}$ is a sequence of appropriate permutation matrices such that, given a vector $\mathbf{x}$, the first elements of the vector $\mathbf{P}^\ite \mathbf{x}$ correspond to the indices in ${I}^k$ and the last elements to the indices in ${A}^k$, for all $k$. By again following a similar reasoning to the one of~\cite[Section 3.3]{Griesse2008}, it can be shown that Line~\ref{eq:ssn_in_algo} of Algorithm~\ref{algo:stackevmfbapp} can be rewritten in such a way that this algorithm is easily seen to be equivalent to an active-set method. In fact, that line is given by 
	\begin{equation*}
	\pr^\iite \leftarrow \inv{\left(\mathbf{P}^\ite\right)} \begin{bmatrix}
	\left([\mathbf{H}]^*_{:{I}^\ite} [\mathbf{H}]_{:{I}^\ite}^{} \right)^{-1} \left[ \mathbf{H}^* \mathbf{b} - \duo^\ite_1 + \tau \mathbf{e}^\ite_{\pm}  \right]_{{I}^\ite}^{} \\
	\mathbf{0}
	\end{bmatrix},
	\end{equation*}
	where $\mathbf{e}^\ite_{\pm} \eqdef \text{sgn} \left[ \pr^k - 2 \tau \left( \mathbf{H}^* \left(\mathbf{H} \mathbf{x}^k - \mathbf{b}\right) + \duo^\ite_1 \right) \right]$, for every $\ite$. The dimension of the problem to solve at each iteration is given by the cardinality of the set $I^\ite$. Naturally, the sparser the solution is estimated to be, the smaller the dimension of this problem is. In contrast, methods such as the \ac{ADMM}~\cite{Afonso2010} require the solution of a problem involving the full matrix $\mathbf{H}^*\mathbf{H}$. This is the reason why semimooth Newton methods are able to achieve faster convergence rates in practice than others.
	
\begin{figure}[!t]
		\begin{center}
			\includegraphics[scale=.45]{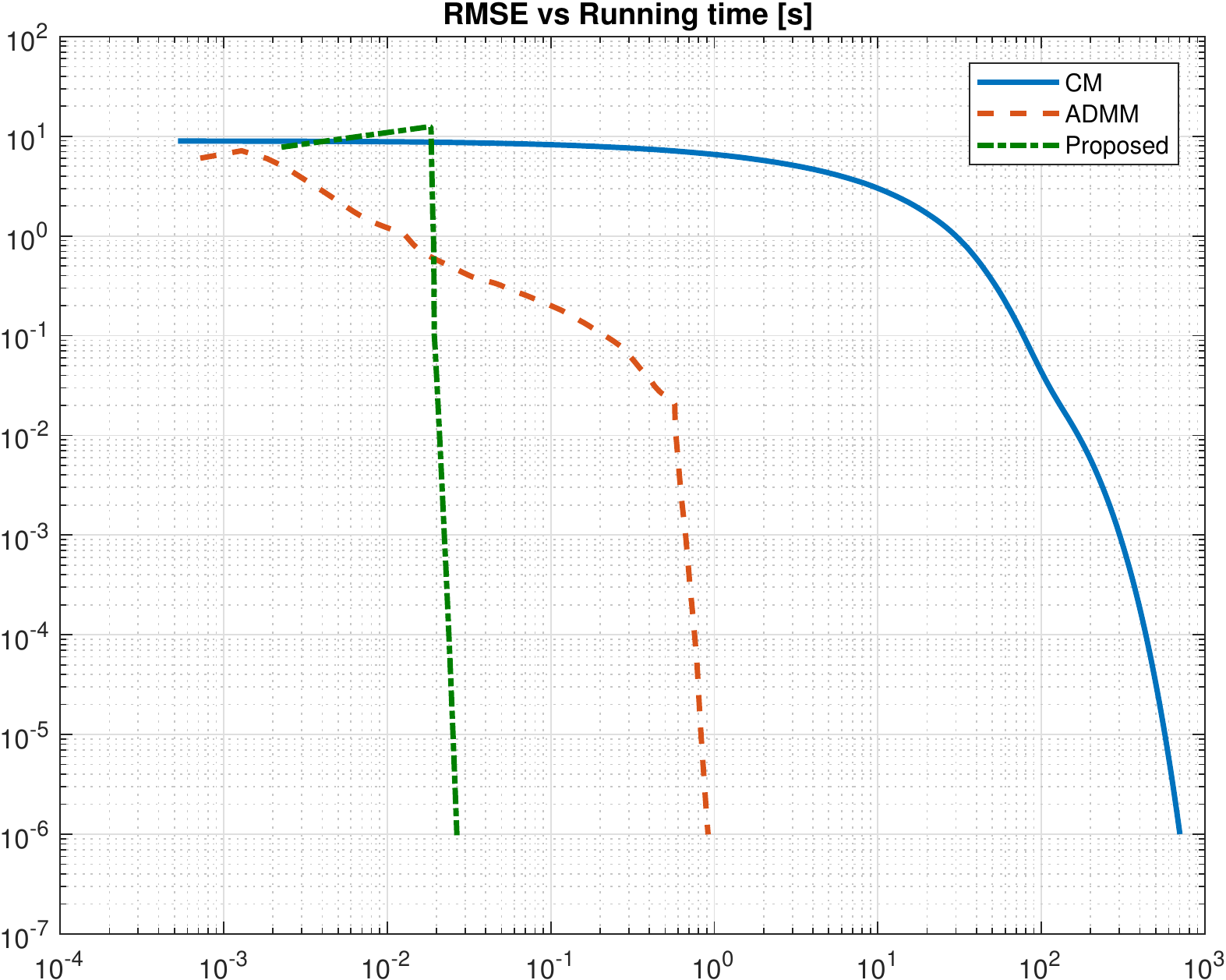}%
			\caption{RMSE, as a function of time, between the estimates of each iteration and the representative solution, for the three methods.}
			\label{fig:rmse_vs_time_inv_int}
		\end{center}	
	\end{figure}

	We simulate an example similar to the one studied by Griesse and Lorenz~\cite[Section 4.1]{Griesse2008} but consider the noise to be Gaussian with a \ac{SNR} of 30 dB. We have set $\mu = 3 \times 10^{-3}$, $c=-80$, and $d=52$. We compared Algorithm~\ref{algo:stackevmfbapp} (denoted in what follows as \emph{Proposed}) with \ac{ADMM} and with the \ac{CM} to solve~\eqref{eq:lasso_constraint}. We manually tuned the different parameters of the three methods in order to achieve the fastest convergence results in practice. We arbitrarily chose the result of \ac{ADMM} after $10^7$ iterations as representative of the solution given by the three methods. Fig.~\ref{fig:rmse_vs_time_inv_int} illustrates the behavior of the three methods by showing the \ac{RMSE} %\footnote{\label{fn:rmse} The \ac{RMSE} is defined as
	%	\begin{equation}
	%	RMSE =  \sqrt{\frac{1}{\dimn} \sum_{i=1}^\dimn \left( [\mathbf{x}]_i - [\mathbf{x}^o]_i \right)^2},
	%	\end{equation}		
	%	where $\mathbf{x}_i$ and $\mathbf{x}_i^o$ denote, respectively, the estimated $\mathbf{x}$ and the original $\mathbf{x}$.} 
	between the estimates of each method and the representative solution, as a function of time. The three methods were initialized with the zero vector. The experiments were performed using MATLAB on an Intel Core i7 CPU running at 3.20~GHz, with 32~GB of RAM. 
%	
%	\begin{figure}[!t]
%		\vspace{-10pt}
%		\centering
%		\subfloat[]{\includegraphics[scale=.4]{u_real.eps}%
%			\label{fig:u_real}}
%		\hfil
%		\subfloat[]{\includegraphics[scale=.4]{observed.eps}%
%			\label{fig:observed}}
%		\hfil
%		\subfloat[]{\includegraphics[scale=.4]{u_estimated.eps}%
%			\label{fig:u_estimated}}	
%		\caption{Original $\mathbf{x}$ (a), observed $\mathbf{b}$ (b), and estimated $\mathbf{x}$ (c).}
%		\label{fig:u}
%		\vspace{-15pt}
%	\end{figure}
	
%	\vspace{-2in}
	
	In this example, we did not enforce assumptions~\eqref{eq:stackevmfbassave} but verified in practice that they were satisfied. However, in more complex examples, it may be necessary to devise a strategy that generates a sequence $\seq{\ave^\ite}$ satisfying these assumptions. This is akin to the necessity of devising globalization strategies in other Newton-like methods~\cite[Chapter 8]{Facchinei2003}.
	
	\subsection{Appraisal}
	
	It is clear that, for this example, the proposed method has a much faster convergence than either \ac{CM} or \ac{ADMM}. This improvement in convergence is similar to the one observed in the methods discussed in Section~\ref{sec:ssn}. In general, the sparser the solution is, the faster the method is as well. In order to benefit from this property, we must be able to solve the lower-dimensional linear system faster than the full system. This may not always be possible: for example, in problems that involve computations with the \ac{FFT} of a signal, we usually have only modest improvements in speed if we wish to compute only selected elements of the \ac{FFT}.\footnote{See~\url{http://www.fftw.org/pruned.html} for details.} However, for large-scale problems and for highly-sparse signals, methods known as sparse \acp{FFT}~\cite{Gilbert2014} may be useful. We verified in other experiments not detailed here that the proposed method has a comparable convergence speed to ADMM in problems whose solutions are not sparse or where we cannot take advantage of their sparsity.
		
	\section{Conclusions} \label{sec:conclusions}
		
	In this work, we defined operator-weighted averaged operators, and showed that they can be used to construct a number of algorithms with good convergence properties. These algorithms have very broad applications, and seem to be particularly suitable to address problems with sparsity-inducing regularizers, as suggested by a simple experiment. Possible future directions to be explored are the possibility of relaxing the assumptions on $\ave^\ite$, and the study of which problems are most suitable to be tackled by these methods.

	\bibliographystyle{IEEEtran}
	% argument is your BibTeX string definitions and bibliography database(s)
	\bibliography{refs}

% Generated by IEEEtran.bst, version: 1.14 (2015/08/26)
\begin{thebibliography}{10}
\providecommand{\url}[1]{#1}
\csname url@samestyle\endcsname
\providecommand{\newblock}{\relax}
\providecommand{\bibinfo}[2]{#2}
\providecommand{\BIBentrySTDinterwordspacing}{\spaceskip=0pt\relax}
\providecommand{\BIBentryALTinterwordstretchfactor}{4}
\providecommand{\BIBentryALTinterwordspacing}{\spaceskip=\fontdimen2\font plus
\BIBentryALTinterwordstretchfactor\fontdimen3\font minus
  \fontdimen4\font\relax}
\providecommand{\BIBforeignlanguage}[2]{{%
\expandafter\ifx\csname l@#1\endcsname\relax
\typeout{** WARNING: IEEEtran.bst: No hyphenation pattern has been}%
\typeout{** loaded for the language `#1'. Using the pattern for}%
\typeout{** the default language instead.}%
\else
\language=\csname l@#1\endcsname
\fi
#2}}
\providecommand{\BIBdecl}{\relax}
\BIBdecl

\bibitem{Figueiredo2003}
M.~Figueiredo and R.~Nowak, ``An {EM} algorithm for wavelet-based image
  restoration,'' \emph{IEEE Trans. Image Process.}, vol.~12, no.~8, pp.
  906--916, Aug 2003.

\bibitem{Daubechies2004}
I.~Daubechies, M.~Defrise, and C.~De~Mol, ``An iterative thresholding algorithm
  for linear inverse problems with a sparsity constraint,'' \emph{Comm. Pure
  Appl. Math.}, vol.~57, no.~11, pp. 1413--1457, 2004.

\bibitem{Moreau1962}
J.~Moreau, ``Fonctions convexes duales et points proximaux dans un espace
  {Hilbertien},'' \emph{Comptes Rendus Acad. Sci.}, vol. A255, pp. 2897--2899,
  1962.

\bibitem{Byrne2004}
C.~Byrne, ``A unified treatment of some iterative algorithms in signal
  processing and image reconstruction,'' \emph{Inverse Probl.}, vol.~20, no.~1,
  pp. 103--120, 2004.

\bibitem{Rockafellar1970}
R.~Rockafellar, \emph{Convex Analysis}.\hskip 1em plus 0.5em minus 0.4em\relax
  New Jersey, USA: Princeton University Press, 1970.

\bibitem{Bauschke2011}
H.~Bauschke and P.~Combettes, \emph{Convex Analysis and Monotone Operator
  Theory in Hilbert Spaces}.\hskip 1em plus 0.5em minus 0.4em\relax New York,
  NY, USA: Springer, 2011.

\bibitem{Schmidt2011}
M.~Schmidt, D.~Kim, and S.~Sra, ``Projected {Newton}-type methods in machine
  learning,'' in \emph{Optimization for Machine Learning}.\hskip 1em plus 0.5em
  minus 0.4em\relax MIT Press, 2011, pp. 305--330.

\bibitem{Becker2012}
S.~Becker and M.~Fadili, ``A quasi-{Newton} proximal splitting method,'' in
  \emph{Proc. 25th Int. Conf. Neural Informat. Process. Systems}, Lake Tahoe,
  Nevada, 2012, pp. 2618--2626.

\bibitem{Lee2014}
J.~Lee, Y.~Sun, and M.~Saunders, ``Proximal {Newton}-type methods for
  minimizing composite functions,'' \emph{SIAM J. Optim.}, vol.~24, no.~3, pp.
  1420--1443, 2014.

\bibitem{Combettes2011b}
P.~Combettes and J.-C. Pesquet, ``Primal--dual splitting algorithm for solving
  inclusions with mixtures of composite, {Lipschitzian}, and parallel-sum type
  monotone operators,'' \emph{Set-Valued Var. Anal.}, vol.~20, no.~2, pp.
  307--330, 2011.

\bibitem{Simoes2017}
\BIBentryALTinterwordspacing
M.~Sim\~{o}es, \emph{On some aspects of inverse problems in image
  processing}.\hskip 1em plus 0.5em minus 0.4em\relax Universidade de Lisboa,
  Instituto Superior T\'{e}cnico, Portugal \& Universit\'{e} Grenoble Alpes,
  France: PhD dissertation, 2017. [Online]. Available:
  \url{http://cascais.lx.it.pt/\%7Emsimoes/dissertation/}
\BIBentrySTDinterwordspacing

\bibitem{Hintermuller2003}
M.~Hinterm\"{u}ller, K.~Ito, and K.~Kunisch, ``The primal--dual active set
  strategy as a semismooth {N}ewton method,'' \emph{SIAM J. Optim.}, vol.~13,
  no.~3, pp. 865--888, 2003.

\bibitem{Facchinei2003}
F.~Facchinei and J.-S. Pang, \emph{Finite-Dimensional Variational Inequalities
  and Complementarity Problems, Vols. I \& II}.\hskip 1em plus 0.5em minus
  0.4em\relax Springer-Verlag, 2003.

\bibitem{Griesse2008}
R.~Griesse and D.~Lorenz, ``A semismooth {N}ewton method for {T}ikhonov
  functionals with sparsity constraints,'' \emph{Inverse Probl.}, vol.~24,
  no.~3, p. 035007, 2008.

\bibitem{Qi1993b}
L.~Qi, ``Convergence analysis of some algorithms for solving nonsmooth
  equations,'' \emph{Math. Oper. Res.}, vol.~18, no.~1, pp. 227--244, 1993.

\bibitem{Ulbrich2011}
M.~Ulbrich, \emph{Semismooth {Newton} Methods for Variational Inequalities and
  Constrained Optimization Problems in Function Spaces}.\hskip 1em plus 0.5em
  minus 0.4em\relax Philadelphia, PA: MOS-SIAM Ser. Optim., 2011.

\bibitem{Combettes2013}
P.~Combettes and B.~V\~{u}, ``Variable metric quasi-{F}ej\'{e}r monotonicity,''
  \emph{Nonlinear Anal-Theor}, vol.~78, pp. 17--31, 2013.

\bibitem{Combettes2014b}
P.~Combettes and B.~V{\~{u}}, ``Variable metric forward--backward splitting
  with applications to monotone inclusions in duality,'' \emph{Optim.},
  vol.~63, no.~9, pp. 1289--1318, 2014.

\bibitem{Afonso2010}
M.~Afonso, J.~Bioucas-Dias, and M.~Figueiredo, ``Fast image recovery using
  variable splitting and constrained optimization,'' \emph{IEEE Trans. Image
  Process.}, vol.~19, no.~9, pp. 2345--2356, Sept 2010.

\bibitem{Condat2013}
L.~Condat, ``A primal--dual splitting method for convex optimization involving
  {L}ipschitzian, proximable and linear composite terms,'' \emph{J. Optim.
  Theory Appl.}, vol. 158, no.~2, pp. 460--479, 2013.

\bibitem{Gilbert2014}
A.~Gilbert, P.~Indyk, M.~Iwen, and L.~Schmidt, ``Recent developments in the
  sparse {Fourier} transform: A compressed {Fourier} transform for big data,''
  \emph{IEEE Signal Process. Mag.}, vol.~31, no.~5, pp. 91--100, Sept 2014.

\end{thebibliography}
	%
	% <OR> manually copy in the resultant .bbl file
	% set second argument of \begin to the number of references
	% (used to reserve space for the reference number labels box)
	
	%\begin{thebibliography}{00}
	%\bibitem{b1} G. Eason, B. Noble, and I. N. Sneddon, ``On certain integrals of Lipschitz-Hankel type involving products of Bessel functions,'' Phil. Trans. Roy. Soc. London, vol. A247, pp. 529--551, April 1955.
	%\bibitem{b2} J. Clerk Maxwell, A Treatise on Electricity and Magnetism, 3rd ed., vol. 2. Oxford: Clarendon, 1892, pp.68--73.
	%\bibitem{b3} I. S. Jacobs and C. P. Bean, ``Fine particles, thin films and exchange anisotropy,'' in Magnetism, vol. III, G. T. Rado and H. Suhl, Eds. New York: Academic, 1963, pp. 271--350.
	%\bibitem{b4} K. Elissa, ``Title of paper if known,'' unpublished.
	%\bibitem{b5} R. Nicole, ``Title of paper with only first word capitalized,'' J. Name Stand. Abbrev., in press.
	%\bibitem{b6} Y. Yorozu, M. Hirano, K. Oka, and Y. Tagawa, ``Electron spectroscopy studies on magneto-optical media and plastic substrate interface,'' IEEE Transl. J. Magn. Japan, vol. 2, pp. 740--741, August 1987 [Digests 9th Annual Conf. Magnetics Japan, p. 301, 1982].
	%\bibitem{b7} M. Young, The Technical Writer's Handbook. Mill Valley, CA: University Science, 1989.
	%\end{thebibliography}
	
	\acrodef{ADMM}{alternating-direction method of multipliers}
	\acrodef{BCCB}{block-circulant-circulant-block}
	\acrodef{BSNR}{blurred-signal-to-noise ratio}
	\acrodef{i.i.d.}{independent and identically distributed}
	\acrodef{KKT}{Karush-Kuhn-Tucker}
	\acrodef{PD}{positive-definite}
	\acrodef{PSD}{positive semidefinite}	
	\acrodef{VMPD}{variable-metric primal--dual method}
	\acrodef{FFT}{fast Fourier Transform}
	\acrodef{SNR}{signal-to-noise ratio}
	\acrodef{RMSE}{root-mean-squared error}
	\acrodef{CM}{method by Condat~\cite{Condat2013}}		
	
\end{document}